\documentclass[12pt,letterpaper]{amsart}

\usepackage{amsmath,amssymb,amsthm}
\usepackage{geometry}
\usepackage{hyperref}

\usepackage{tikz}
\usetikzlibrary{arrows.meta}

\usepackage[all]{xy}

\usepackage{graphicx}

\usepackage{float, tabularx, booktabs, caption}

\usepackage{epsfig}

\newtheorem{innercustomgeneric}{\customgenericname}
\providecommand{\customgenericname}{}
\newcommand{\newcustomtheorem}[2]{%
  \newenvironment{#1}[1]
  {%
   \renewcommand\customgenericname{#2}%
   \renewcommand\theinnercustomgeneric{##1}%
   \innercustomgeneric
  }
  {\endinnercustomgeneric}
}

\newcustomtheorem{customthm}{Theorem}

\newtheorem*{theorem}{Theorem}

\newtheorem{definition}{Definition}

\numberwithin{equation}{section}

\title[Robust transitivity without sectional-hyperbolicity]{Robust transitivity without sectional-hyperbolicity}

\author{A. Arbieto}
\address{Instituto de Matem\'atica, Universidade Federal do Rio de Janeiro, Rio de Janeiro, Brazil.}
\email{arbieto@im.ufrj.br}

\author{W. Britto}
\address{Instituto de Matem\'atica, Universidade Federal do Rio de Janeiro, Rio de Janeiro, Brazil.}
\email{walter\_britto@im.ufrj.br}

\author{C.A. Morales}
\address{Hangzhou International Innovation Institute of Beihang University, Hangzhou 311115, China.}
\email{morales@impa.br}

\author{E. Rego}
\address{Faculty of Applied Mathematics, AGH University of Science and Technology, Krakow, Poland.}
\email{rego@agh.edu.cn}


\date{}

\begin{document}

\begin{abstract}
For any integer $n\geq 5$, we construct an $n$–dimensional $C^1$ vector field exhibiting a robustly transitive singular attractor which is not sectional–hyperbolic. Nevertheless, the attractor is singular–hyperbolic. This provides the first such examples improving some features of the constructions in \cite{lyyz,ts}.
\end{abstract}

\maketitle

\section{Introduction}

\noindent
It is well known that every robustly transitive attractor of a $3$–dimensional $C^1$ vector field is sectional–hyperbolic \cite{mpp1}. We believe that this phenomenon is specific to dimension three and does not extend to higher dimensions. More precisely, we expect that in dimension $4$ (and hence in every dimension $\geq 4$) there exist $C^1$ vector fields exhibiting robustly transitive attractors that are not sectional–hyperbolic. Of course, this question is meaningful only in the singular setting (by Theorem~B in \cite{Ma}). Here, singular means that the set contains singularities, all of them hyperbolic.
If in the above question we replace ``robustly transitive attractor'' with ``robustly chain transitive attractor'', then the answer is positive by theorems A and B in \cite{bly}.

Partial evidence in both directions is available. On the one hand, Theorem~C in \cite{sgw} shows that every singular attractor of a $4$–dimensional $C^1$ vector field far from nonhyperbolic closed orbits is sectional–hyperbolic, thus pointing toward a negative answer. Likewise,
by \cite{mm1} (see also Theorem~A in \cite{mm}), in any dimension, a robustly transitive singular attractor which is strongly homogeneous of index $d$
and whose singularities have index $>d$ must be sectional–hyperbolic
(this last condition was removed in \cite{zgw}).

On the other hand, two remarkable examples of robustly Lyapunov stable chain recurrent classes associated with hyperbolic singularities of certain $4$–dimensional vector fields suggest the opposite conclusion. The first one, constructed in \cite{ts} and known as the {\em wild strange attractor}, is singular–hyperbolic but not sectional–hyperbolic due to the presence of homoclinic tangencies. The second example, recently obtained in \cite{lyyz} and called a {\em derived–from–Lorenz} attractor, is far from homoclinic tangencies but fails to be sectional–hyperbolic since it robustly contains a heterodimensional cycle.

For convenience, let us recall the relevant terminology. A compact invariant set $\Lambda$ of a $C^1$ vector field $X$ is said to be {\em strongly homogeneous of index $d$} \cite{lgw} if there exist a $C^1$ neighborhood $\mathcal{U}$ of $X$ and a neighborhood $U$ of $\Lambda$ such that the index of the periodic orbits contained in $U$ is $d$ for all vector fields in $\mathcal{U}$.

On the other hand, we say that the chain recurrent class $C$ associated with a hyperbolic closed orbit $O$ of $X$ is {\em robustly Lyapunov stable} if there exist a neighborhood $U$ of $C$ and a $C^1$ neighborhood $\mathcal{U}$ of $X$ such that, for every $Y \in \mathcal{U}$, the chain recurrent class $C_Y$ of the continuation $O_Y$ is the unique Lyapunov stable set of $Y$ contained in $U$.

Likewise, $C$ is said to {\em robustly contain a heterodimensional cycle} if there exists a $C^1$ neighborhood $\mathcal{U}$ of $X$ such that, for every $Y \in \mathcal{U}$, the chain recurrent class $C_Y$ of the continuation $O_Y$ contains two transitive hyperbolic sets $\Lambda_1$ and $\Lambda_2$ with different indices whose continuations satisfy
\[
W^u(\Lambda_{1,Z}) \cap W^s(\Lambda_{2,Z}) \neq \emptyset
\quad \text{and} \quad
W^u(\Lambda_{2,Z}) \cap W^s(\Lambda_{1,Z}) \neq \emptyset
\]
for every vector field $Z$ sufficiently close to $Y$.

Although both examples mentioned above contain singularities and fail to be sectional–hyperbolic, it remains unknown whether they are attracting, transitive, or robustly transitive (in the case where they are attracting). Consequently, the existence of robustly transitive singular attractors that are not sectional–hyperbolic remains open, at least in dimension $\geq4$.

The purpose of this paper is to provide a positive answer to this question in dimension $\geq 5$. More precisely, for every $n \geq 5$, we construct an $n$–dimensional $C^1$ vector field exhibiting a robustly transitive singular attractor that is not sectional–hyperbolic. Nevertheless, we prove that this attractor is {\em singular–hyperbolic}, a notion introduced in \cite{m2,mpp,ts}, which is sometimes confused in the literature with sectional–hyperbolicity (see, for instance, \cite{bdl,bly,CWYZ,clyz,wwy,zgw}).\footnote{Sectional–hyperbolicity was originally introduced in the preprint \cite{mm1} which was published later in \cite{mm}. It is strictly stronger than singular–hyperbolicity  (see also \cite{Sal}).}
Actually, it seems that every robustly transitive singular attractor of a
\(C^1\) vector field is singular-hyperbolic.

Let us present our results in a precise way.
\noindent
By an {\em $n$-dimensional $C^1$ vector field} we mean a vector field $X:M\to TM$ of class $C^1$ defined on a closed (i.e. compact connected boundaryless Riemannian) manifold $M$ of dimension $n$. Its flow is denoted by $X_t$.
We say that $\sigma\in M$ is a singularity if $X_t(\sigma)=\sigma$ for all $t\in\mathbb{R}$. We say that $x\in M$ is a periodic orbit is there is a minimal positive time $t$ (called period) such that $X_t(x)=x$.
The {\em orbit} of $x\in M$ is defined by $O(x)=\{X_t(x)\mid t\in\mathbb{R}\}$.
The orbit of a periodic point will be referred to as a periodic orbit.
Notice that an orbit (i.e. the orbit of some point) is closed (hence compact) if and only if it is either a singularity or a periodic orbit. We denote by $Sing(X)$ the set of singularities of $X$.

We say that $\Lambda\subset M$ is {\em invariant} if $X_t(\Lambda)=\Lambda$ for all $t\in\mathbb{R}$.
A compact invariant set is {\em transitive} if for every pair of open sets $U,V$ intersecting $\Lambda$ and $T>0$ there is $t\geq T$ such that
$X_t(U)\cap V\cap \Lambda\neq\emptyset$.
On the other hand, $\Lambda$ is {\em dynamically isolated} if it exhibits a compact neighborhood $U$ (called isolating block) such that
$$
\Lambda=\bigcap_{t\in\mathbb{R}}X_t(U).
$$
A dynamically isolated set $\Lambda$ is {\em robustly transitive} if
there is an isolating block $U$ such that
$\Lambda_Y$ is transitive for every $C^1$ vector field $Y$ which is $C^1$ close to $X$ where
$$
\Lambda_Y=\bigcap_{t\in\mathbb{R}}Y_t(U)
$$
is the continuation of $\Lambda$ for $Y$ close to $X$.
A dynamically isolated set is {\em attracting} if it exhibits a positively invariant isolating block $U$, namely,
$X_t(U)\subset U$ for all $t\geq0$.
An {\em attractor} is a transitive attracting set.

\begin{definition}
A compact invariant set $\Lambda$ is {\em hyperbolic} if there is a tangent bundle decomposition $T_\Lambda M=E^s_\Lambda\oplus E^u_\Lambda\oplus E^u_\Lambda$ (called {\em hyperbolic splitting}) along with positive constants $C,\lambda$ such that
$E^X_x=\langle X(x)\rangle$ is the subspace generated by $X(x)$, $E^*_\Lambda$ is invariant 
(i.e. $DX_t(x)E^*_x=E^*_{X_t(x)}$ for $*=s,u$ and $x\in \Lambda$),
$$
\|DX_t(x)|_{E^s_x}\|\leq Ce^{-\lambda t}\quad\mbox{ and }\quad
m(DX_t(x)|_{E^u_x})\geq C^{-1}e^{\lambda t},\qquad\forall x\in\Lambda,\, t\geq0,
$$
where $m(\cdot)$ stands for the minimal norm of a linear operator (also called mininorm).
\end{definition}

A closed orbit $O$ is hyperbolic if it does as a compact invariant set.
A compact invariant set is {\em singular} if it has singularities, all of them hyperbolic.

The invariant manifold theory \cite{HPS} asserts that through any point $x$ of a hyperbolic set it passes a {\em stable manifold} $W^{s}(x)$ and a {\em unstable manifold} $W^{u}(x)$ tangent at $x$ to the stable and unstable directions $E^s_x$ and $E^u_x$ respectively. The stable manifold along the orbit $O(x)=\{X_t(x)\mid t\in\mathbb{R}\}$ is defined by
$$
W^u(O(x))=\bigcup_{t\in\mathbb{R}}W^u(X_t(x)).
$$

An invariant tangent bundle decomposition $T_\Lambda M=E_\Lambda\oplus F_\Lambda$ over a compact invariant set $\Lambda$ is {\em dominated} if there are positive constants $C,\lambda$ such that
$$
\frac{\|DX_t(x)|_{E_x}\|}{m(DX_t(x)|_{F_x})}\leq Ce^{-\lambda t},\qquad\forall x\in\Lambda,\, t\geq0.
$$
We say that $\Lambda$ is a {\em partially hyperbolic set} if it exhibits a partially hyperbolic splitting. The latter is a dominated splitting $T_\Lambda M=E^s_\Lambda\oplus E^c_\Lambda$ with contracting dominating direction namely $\exists C,\lambda>0$ such that
$$
\|DX_t(x)|_{E^s_x}\|\leq Ce^{-\lambda t},\qquad\forall x\in \Lambda,\, t\geq0.
$$

\begin{definition}
\label{cucu}
A compact invariant set $\Lambda$ is \emph{singular-hyperbolic} if its singularities are hyperbolic and there exists a partially
hyperbolic splitting $T_\Lambda M=E^s_\Lambda\oplus E^c_\Lambda$ for the tangent flow $DX_t$ such that the
central direction is \emph{volume expanding}, namely, there exist $C,\lambda>0$
such that
\[
\big|\det\big(DX_t(x)\big|_{E^c_x}\big)\big|
\ge Ce^{\lambda t},
\qquad\forall x\in\Lambda,\ t\ge 0.
\]
\end{definition}

Every hyperbolic set is singular-hyperbolic. The basic example of a singular-hyperbolic set which is not hyperbolic is the {\em geometric Lorenz attractor} \cite{abs,lo,g}.

\begin{definition}
A compact invariant set $\Lambda$ is {\em sectional-hyperbolic} if its singularities are hyperbolic and there is a partially hyperbolic splitting $T_\Lambda M=E^s_\Lambda\oplus E^c_\Lambda$ with {\em sectionally expanding central direction}, namely, there exist $C,\lambda>0$ such that
$$
|det(DX_t(x)|_{L_x}|\geq C^{-1}e^{\lambda t},\qquad\forall t\geq0,
$$
for all $x\in X$ and all two-dimensional subspace $L_x\subset E^c_x$.
\end{definition}

Apart from the geometric Lorenz attractors, the basic examples of a sectional-hyperbolic sets are the {\em multidimensional Lorenz attractors} \cite{BPV}.
Any sectional-hyperbolic set is singular-hyperbolic but not conversely as, for instance, the wild strange attractor.
With these definitions we can state our result.

\begin{theorem}
For every integer $n\geq5$ there is an $n$-dimensional $C^1$ vector field exhibiting a robustly transitive singular attractor which is not sectional-hyperbolic. Moreover, this attractor is singular-hyperbolic.
\end{theorem}

Let us comment on the proof of this theorem.
First, observe that it suffices to establish the result in dimension \( n = 5 \), since examples in higher dimensions can be obtained by taking the product with suitable contracting directions. To prove it for $n=5$
we will insert a singularity, following \cite{BPV} or \cite{n}, into the suspension of a suitable $4$–dimensional solenoid,
rather than perform a DA bifurcation \cite{W} on a $4$–dimensional singular attractor (as the $4$-dimensional Lorenz attractor in \cite{lyyz}).

This strategy requires a careful choice of the underlying solenoid.
We discard the ones with foliation map given by \cite{ar,hg,s}
as the robust transitivity in these cases rely on invariant subbundles that maybe destroyed after inserting the singularity.

Instead, we proceed as in Example~7 of \cite{ass} starting with an endomorphism \( f : T \to T \) of the two–dimensional torus \( T \). In contrast with \cite{ass}, where \(f\) is chosen to be expanding, we take \(f\) as in Section~5 of \cite{lp}, assuming volume expansion among another properties (another possible choice is given in \cite{r}). The idea is to apply the results of \cite{lp}, which provide general conditions ensuring that a volume expanding endomorphism of \(T\) is robustly transitive.

The desired solenoid is then given by an embedding
\[
F : T \times D \longrightarrow T \times D ,
\]
which contracts the vertical foliation \( * \times D \), while the dynamics along the base are governed by \(f\).

Once this solenoid is chosen we suspend it to obtain a $C^1$ vector field on a (5)–dimensional manifold and insert the singularity as already mentioned.

The foliation map \( f^* \) obtained after the insertion of the singularity satisfies the hypotheses for robust transitivity in \cite{lp}, except that it is defined on a punctured torus. Nevertheless, the methods of \cite{lp} still apply and imply that \( f^* \) is robustly transitive. This yields robust transitivity of the corresponding attracting set and ultimately produces the desired vector field.

\section{Proof of the theorem}

\noindent
As already said it suffices to prove the theorem for $n=5$.
The proof is organized in several steps.

\subsection{Choice of toral endomorphisms}

We say that an endomorphism $f:T\to T$ satisfies {\rm (LP)} if the following properties hold:

\begin{enumerate}
\item[(LP1)]
$f$ is volume expanding, i.e., there exists $\Delta>1$ such that
\[
|\det Df(x)|>\Delta \quad \text{for all } x\in T.
\]

\item[(LP2)]
For every $x\in T$, its preorbit
\[
\bigcup_{n\ge0} f^{-n}(x)
\]
is dense in $T$.

\item[(LP3)]
There exists an open set $U_0\subset T$ with $\operatorname{diam}(U_0)<1$ such that
$f$ is uniformly expanding on its complement $U_0^c=T\setminus U_0$, namely, there exists $\Delta>1$ such that
\[
\|Df(x)v\|\ge \Delta\|v\|
\]
for every $x\in U_0^c$ and every $v\in T_xT$.

\item[(LP4)]
There exist $\delta_0>0$ small and an open neighborhood $U_1$ of $\overline{U_0}$ (the closure of $U_0$) such that every arc $\gamma\subset U_0^c$ with diameter greater than $\delta_0$ contains a point $y$ satisfying
\[
f^k(y)\in U_1^c \quad \text{for all } k\ge1.
\]

\item[(LP5)]
$U_1^c$ is forward invariant, i.e.,
\[
U_1^c\subset f(U_1^c).
\]
\end{enumerate}

It was proved in \cite{lp} that every endomorphism satisfying {\rm (LP)} is $C^1$–robustly transitive.

Fix such an endomorphism $f:T\to T$. We further assume:

\begin{enumerate}
\item[(H6)]
$f$ has two fixed points, a repelling one in $U_1^c$ and a saddle in $U_0$.
\end{enumerate}

Examples of such endomorphisms can be constructed as in Section~5 of \cite{lp}.

\subsection{The solenoid}

Consider a $C^1$ embedding
\[
F:T\times D\to T\times D
\]
preserving the vertical foliation $\{x\}\times D$ and covering the chosen endomorphism $f$. That is,
\[
\pi_1\circ F = f\circ \pi_1,
\]
where $\pi_1:T\times D\to T$ is the projection.
Assume moreover that $F$ contracts uniformly the vertical foliation.

Then the vertical bundle
\[
E^v=\ker D\pi_1
\]
is a continuous $DF$–invariant contracting subbundle over the maximal invariant set
\[
\Lambda=\bigcap_{n\ge0} F^n(T\times D).
\]
Since $E^v$ has constant dimension and is uniformly contracted, by choosing the vertical contraction sufficiently strong compared with the transverse dynamics induced by $f$ and using graph-transformed techniques \cite{HPS}, there exists a continuous complementary bundle $E^c$ such that
\begin{equation}
\label{splitF}
T_\Lambda(T\times D)
=
E^v\oplus E^c
\end{equation}
is a dominated splitting.
Furthermore, since $f$ is volume expanding by (LP1), the bundle $E^c$ is volume expanding too.

\subsection{Suspension}

Let $M$ be the quotient manifold obtained from
$T\times D\times[0,1]$
by identifying
\[
(z,1)\sim (F(z),0).
\]
Consider the suspension vector field
\[
\hat X=\frac{\partial}{\partial t}.
\]
The cross section
\[
\Sigma=T\times D\times\{0\}
\]
has return map $F$.

Transporting the splitting \eqref{splitF} along the flow yields an invariant splitting
\begin{equation}
\label{splitSusp}
T_{\hat\Lambda}M
=
\hat E^v
\oplus
E^{\hat X}
\oplus
\hat E^c,
\end{equation}
over
\[
\hat\Lambda
=
\bigcap_{t\ge0}\hat X_t(M),
\]
where
\[
\dim(\hat E^v)=\dim(\hat E^c)=2,
\]
$\hat E^v$ is uniformly contracting, and $\hat E^c$ is volume expanding.

\subsection{Insertion of a singularity}
Choose a point $q\in U_0$ (see (LP2)), consider a small neighborhood $U\subset U_0$ of $q$ and let
\[
V=U\times D\times[0,1].
\]
Modifying $\hat{X}$ inside $V$, as in the Cherry–flow construction \cite{pdm}, we create a vector field $X$ exhibiting a hyperbolic singularity $\sigma$ with two expanding eigenvalues and three contracting eigenvalues.
This modification is performed so that:

\begin{itemize}
\item the stable manifold of $\sigma$ intersects $\Sigma$ transversely along
\[
\Gamma=\{q\}\times D\times\{0\},
\]
\item the unstable manifold returns to $\Sigma$ along a circle $S_\sigma\subset U$.
\end{itemize}
Define
\[
A
=
\bigcap_{T\ge0}
\overline{\bigcup_{t\ge T} X_t(\Sigma)}.
\]
It follows that $A$ is a connected attracting set of $X$.
In the sequel, we prove that this set is singular-hyperbolic, robustly transitive but not sectional-hyperbolic.

\subsection{Singular- but not sectional-hyperbolicity}

As already said, $A$ is a connected attracting set of $X$.
If it were sectional-hyperbolic, then it could not contain periodic orbits of different indices. As it contains such periodic orbits by assumption {\rm (H6)}, we conclude that $A$ cannot be sectional-hyperbolic.

On the other hand, the deformation transforming $\hat X$ into $X$ preserves the invariant splitting \eqref{splitSusp} outside a small neighborhood of the inserted singularity and extends continuously over the whole attracting set. More precisely, it yields an $X$–invariant splitting
\[
T_AM
=
E^s_A
\oplus
E^c_A,
\]
where $E^s_A$ is the continuation of the contracting bundle $\hat E^v$ and $E^c_A$ is the continuation of the volume expanding bundle $E^{\hat X}\oplus \hat E^c$. By construction, $E^s_A$ is uniformly contracting and $E^c_A$ is volume expanding. Therefore, $A$ is singular-hyperbolic.

\subsection{Robust transitivity}

It remains to prove that $A$ is robustly transitive. The proof is based on the analysis of the return map
\[
F^*:\Sigma\setminus\Gamma\to\Sigma
\]
induced by $X$.
It still preserves and strongly contracts the vertical foliation inducing
a base map
\[
f^*:T\setminus\{q\}\to T
\]
which satisfies {\rm (LP)} except that it is not defined in $q$. Notice that $f^*$ is a {\em local diffeomorphism} in the sense that there is a finite open covering of $T^2$ such that 
$f^*$ restricted each element of this covering (minus $\{q\}$) is a diffeomorphism onto its image.

Next we use the main theorem of \cite{lp}, suitably adapted to this $f^*$: by {\rm (LP)} it follows that $f^*$ satisfies the internal radius growth property:

\begin{enumerate}
\item[(IRG)]
There exists $R_0>0$ such that for every open set $U\subset T\setminus\{q\}$ there exist $x\in U$ and $K\ge0$ such that $(f^*)^K(U)$ contains the ball of radius $R_0$ centered at $(f^*)^K(x)$.
\end{enumerate}

Indeed, define
\[
\Lambda_1
=
\bigcap_{n\ge0}(f^*)^n(U_1^c).
\]
Since $z_0\in U_0$, we can assume as in \cite{lp} that it is a dynamically isolated expanding set. Let $R_0$ be the distance between $\Lambda_1$ and $U_0$.

Since $f^*$ is volume expanding by {\rm (LP1)}, its lift to the universal cover $\mathbb{R}^2$ is a volume expanding map which is defined in $\mathbb{R}^2$ minus countably many points over $q$. It follows that the diameter of any open set grows under iteration. Hence, there exists $N\ge0$ such that $(f^*)^N(U)$ intersects $U_0^c$ in a set of diameter larger than $\delta_0$. Therefore, this image contains an arc satisfying the hypothesis of {\rm (LP4)}, and so some forward iterate enters the expanding region $U_1^c$. This implies that the internal radius grows under iteration, proving {\rm (IRG)}.

Now, for a $C^1$ vector field $Z$ close to $X$, the contracting vertical foliation need not be differentiable, in which case the foliation map $g$ induced by $Z$ may not be differentiable. However, this technical issue may be overcame by dealing directly with the return map (like \cite{ccs} or the triangular maps techniques in \cite{bm,bm1}) thus we will not be concerned about this issue in what follows. In particular, the expanding region $U_1^c$ is not affected by the singularity, and therefore the expanding set $\Lambda_1$ admits a continuation $\Lambda_g$ for $g$. As in \cite{lp}, this continuation intersects every sufficiently large arc in $U_1^c$. Using this property and following the arguments of Subsections 2.4–2.6 of \cite{lp}, we conclude that $g$ also satisfies {\rm (IRG)}.

By Subsection 2.7 of \cite{lp}, property {\rm (IRG)} and the fact that $f^*$ is a local diffeomorphism imply the density of the preorbits of $g$ hence $g$ is transitive. Finally, since the vertical foliation persists under perturbations, the transitivity of the return map implies the transitivity of the continuation of $A$ for $Z$ (see \cite{b}). Therefore, $A$ is robustly transitive completing the proof.
\qed

\section*{Funding}

\noindent
A.A. was partially supported by FAPERJ (Cientista do nosso estado E-26/201.181/2022), CAPES, CNPq (Projeto Universal) and PRONEX-Dynamical Systems.
W.B was partially supported by CNPq-Brazil.

\begin{table}[h]
\begin{tabularx}{\linewidth}{p{1.5cm}  X}
\includegraphics [width=1.8cm]{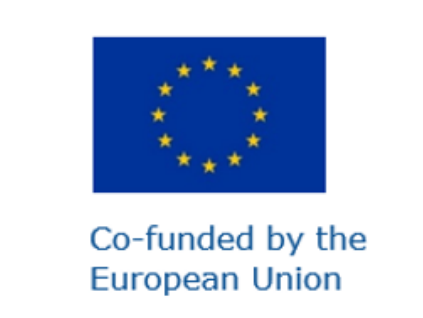} &
\vspace{-1.5cm}
This research is part of a project that has received funding from
the European Union's European Research Council Marie Sklodowska-Curie Project No. 101151716 -- TMSHADS -- HORIZON--MSCA--2023--PF--01.\\
\end{tabularx}
\end{table}

\section{Declaration of competing interest}

\noindent
There is no competing interest.

\section*{Data availability}

\noindent
No data was used for the research described in the article.

\end{document}